\documentclass[11pt, a4paper]{article}
\usepackage{authblk}
\usepackage[utf8]{inputenc}
\usepackage{amsmath}
\usepackage{amssymb}
\usepackage{mathtools}
\usepackage{amsthm}
\usepackage{authblk}
\usepackage{amssymb,enumitem,amsmath, mathtools}
\newtheorem{thm}{Theorem}

\begin{document}

%\begin{frontmatter}



\title{Remarks on the preservation and breaking of translational symmetry for a class of ODEs}

\author[a]{Edward Huynh}
%\ead{huynhe3@unlv.nevada.edu}
\affil{Department of Mathematical Sciences, University of Nevada, Las Vegas, Las Vegas, Box 454020, NV, USA}

%\affiliation[a]{organization={Department of Mathematical Sciences},%Department and Organization
%            addressline={University of Nevada, Las Vegas}, 
%            city={Las Vegas},
%            postcode={Box 454020}, 
%            state={NV},
%            country={USA}}

\author[a]{Keoni Castellano}
%\ead{castek1@unlv.nevada.edu}
\affil{Department of Mathematical Sciences, University of Nevada, Las Vegas, Las Vegas, Box 454020, NV, USA}
\maketitle

\begin{abstract}
In this paper, we provide both a preservation and breaking of symmetry theorem for $2\pi$-periodic problems of the form
\begin{align*}
    \begin{cases}
    -u''(t) + g(u(t)) = f(t) \\
    u(0) - u(2\pi) = u'(0) - u'(2\pi) = 0
    \end{cases}
\end{align*}
where $g: \mathbb{R} \to \mathbb{R}$ is a given $C^1$ function and $f: [0,2\pi] \to \mathbb{R}$ is continuous. We provide a preservation of symmetry result that is analogous to one given by Willem (Willem, 1989) and a generalization of the theorem given by Costa-Fang (Costa and Fang, 2019). Both of these theorems use group actions that are not normally considered in the literature.
\end{abstract}



%\begin{keyword}
%Periodic solutions \sep Symmetry breaking \sep Symmetry preservation %\sep Group action \sep Critical groups \sep Morse index
%\end{keyword}

%\end{frontmatter}


\section{Preliminaries}
\label{sec:prelim}
Let $H$ be a Hilbert space, $(G,\cdot)$ a topological group, and $\{T(g)\mid g\in G\}$ be an isometric representation of $G$ on $H$. In other words, $T(g): H \rightarrow H$ is an isometry satisfying
\begin{enumerate}[label = (\roman*)]
    \item $T(e) = \text{Identity}$
    \item $T(g_1\cdot g_2) = T(g_1) \circ T(g_2)$ for all $g_1,g_2\in G$
    \item $G\times H \ni (g,u) \mapsto T(g)u \in H$ is a continuous map.
\end{enumerate}

The $\textit{orbit}$ of an element $u\in H$ is the set $\{T(g)u\mid g\in G\}$. We say that a subset $A\subset H$ is $\textit{invariant}$ (under $G$) if $T(g)A \subset A$ for all $g\in G$. A functional $\varphi: H \rightarrow \mathbb{R}$ is said to be $\textit{invariant}$ if $\varphi\circ T(g) = \varphi$ for all $g\in G$. A mapping $R: A_1 \rightarrow A_2$ between subsets of $H$ is called $\textit{equivariant}$ if $R\circ T(g) = T(g) \circ R$ for all $g\in G$. 

The set of \textit{fixed points} of $H$ under the representation $\{T(g)\}$, which we refer to as the set of \textit{(most) symmetric elements} of $H$ (under $G$), is the closed subspace of $H$ defined by
\begin{equation*}
    \text{Fix}(G) = \{u\in H\mid T(g)u = u, \ \forall g\in G\}.
\end{equation*}

If $u_0\in H$ is an isolated critical point of $\varphi \in C^2(H, \mathbb{R})$ which is {\it nondegenerate}, i.e. the bilinear form $\varphi''(u_0): H\times H \rightarrow \mathbb{R}$ is nondegenerate, the \textit{Morse Index} of $\varphi$ at $u_0$, which we denote by
\begin{equation*}
    i_M(\varphi,u_0),
\end{equation*}
is the supremum of $k\in \mathbb{N}$ such that $\varphi''(u_0)$ is \textit{negative} {\it definite} on a $k$-dimensional subspace of $H$. 
%The \textit{critical groups} (over a field $\mathbb{F}$) of %$\varphi$ at $u_0$ are defined by
%\begin{equation*}
%    C_n(\varphi,u_0) \coloneqq H_n(\varphi^c\cap \Bar{U}, %\varphi^c\cap \Bar{U}\setminus \{u_0\}), \ n\in \mathbb{N}
%\end{equation*}
%Where $U\subset H$ is an isolating neighborhood of $u_0$, $c = \varphi(u_0),\varphi^c=\{u\in H\mid \varphi(u)\leq c\}$, and $H_n(A,B)$ denotes the $n$th \textit{singular homology group} of the pair $(A,B)$ over the field $\mathbb{F}$. One can show that the critical groups are independent of $U$ and, when $u_0$ is a nondegenerate critical point of $\varphi$ with Morse index $k$, then we have that
%\begin{gather*}
%    C_n(\varphi,u_0) = \delta_{n,k}\mathbb{F}.
%\end{gather*}

Moreover, recall that the eigenvalues of the problem
\begin{gather}\label{eigprob}
    \begin{cases}
    -h''(t) = \lambda h(t)\\
    h(0) = h(2\pi)\\
    h'(0) = h'(2\pi)
    \end{cases}
\end{gather}
are $\lambda_j = j^2$ where $j\in \mathbb{Z}^{0+} = \{z\in \mathbb{Z}\mid z \geq 0\}$, with corresponding eigenfunctions $h_j(t) = C_j\cos(jt) + D_j\sin(jt)$ (and so $h_0(t) = C_0)$, where $C_j,D_j$ are arbitrary constants. Note that $\lambda_0 = 0$ is simple while $\lambda_j, \ j\geq 1$ are double. Furthermore, for $j \geq 1$, the eigenfunctions are periodic with minimal period $\frac{2\pi}{j}$. 

Consider the following spaces
\begin{align*}
    H^1_{per}[0,2\pi] &= \{u\in H^1[0,2\pi]\mid u(0) = u(2\pi), \ u'(0) = u'(2\pi)\};\\
    V_j &= \{u\in H^1_{per}[0,2\pi]\mid \text{u is $2\pi$/j periodic}\};\\
    E_j &= \{C_1\cos(jt) + C_2\sin(jt)\mid C_1,C_2\in \mathbb{R}\}, \ j\in \mathbb{Z}^+;\\
    E_0 &= \mathbb{R}.
\end{align*}
Note that $E_j$ is the $j$th eigenspace of \eqref{eigprob}. We endow $H^1_{per}[0,2\pi]$ with the norm and inner product, respectively,
\begin{gather*}
    ||u||=\left(\int_0^{2\pi}(|u'|^2 + |u|^2)\ dt\right)^{1/2} \ \ \text{and} \ \ \langle u,v \rangle = \int_0^{2\pi} (u'v' + uv)\ dt.
\end{gather*}

\section{Introduction}
\label{sec:intro}
We are interested in problems of the form
\begin{equation}\label{eq:problem}
    \begin{cases}
    -u''(t) + g(u(t)) = f(t)\\
    u(0) = u(2\pi)\\
    u'(0) = u'(2\pi) 
    \end{cases}
\end{equation}
where $g:\mathbb{R} \rightarrow \mathbb{R}$ is a given $C^1$ function and $f:[0,2\pi]\rightarrow \mathbb{R}$ is continuous. Let $G$ be a finite group such that $G = \mathbb{Z}_m$ for some $m$. Suppose $u \in H_{per}^1[0,2\pi]$, then let $\{T(g): g\in G\}$ be an isometric topological representation of $G$ on $H_{per}^1[0, 2\pi]$ such that the action of $T(g)$ on $H_{per}^1[0, 2\pi]$ is given by $T(g)u = \hat{u}(t+ \frac{2\pi g}{m})$ where $\hat{u}$ is the periodic extension of $u$ on the real line. Here 
\begin{align*}
    \text{Fix}(G) = \{u\in H_{per}^1[0,2\pi]\mid u\text{ is $2\pi/m$-periodic}\}.
\end{align*}
For this group action, we call elements of Fix($g$) \textit{translationally symmetric}, since the action preserves symmetry under translation. For such a problem as \eqref{eq:problem}, three questions are of use to us.

The first question concerns that of \textit{preservation of symmetry}. That is, if $f(t)$ has a certain symmetry, say periodicity on the real line with period $P$, then if all solutions to \eqref{eq:problem} also have period $P$, we say that the symmetry is preserved.

Conversely, the second question concerns the \textit{breaking of symmetry}. That is, if $f(t)$ is periodic on the real line with period $P$, then a solution $v(t)$ to $\eqref{eq:problem}$ breaks symmetry if $v(t)$ is not periodic with period $P$.

The third question concerns \textit{existence and multiplicity of solutions using symmetry}. In other words, given an isometric representation of a topological group $G$ on $H$, there exist either one or multiple solutions that are geometrically distinct, i.e. solutions that are not within the orbit of the other.

Results on preservation and breaking of symmetry have appeared for both ordinary and partial differential equations such as in Willem \cite{Willem1989}, Dancer \cite{Dancer1983}, Lazer-McKenna \cite{Lazer1988}, just to mention a few pioneering references on the subject. Costa-Fang \cite{CostaFang2019} proved a result for breaking of symmetry that concerns translational symmetry using $\mathbb{Z}_p$ ($p$ prime) group actions, i.e. they showed that there exists a function that is $2\pi/\hat{p}$ ($\hat{p} \neq p$ and $\hat{p}$ being prime) periodic such that (\ref{eq:problem}) has a solution that does not have the same periodicity. The theorem relies on the hypothesis that the derivative on $g$ must be bounded above and below by certain eigenvalues. In this way, this is a type of "symmetry" imposed on $g$. Our breaking of symmetry theorem extends the hypothesis to cases when the order of the group action and the order of the eigenvalues are relatively prime. On the other hand, our preservation of symmetry result shows that all solutions of (\ref{eq:problem}) must satisfy a certain symmetry if the nonlinear term satisfies a symmetry based on the eigenvalues.

There have also been results that exhibit existence of solutions for Hamiltonian systems that preserve a type of symmetry. Major approaches in this field include Lusternik-Schnirelmann Theory and Index Theory as in Costa-Willem \cite{Costa1986}. In particular, existence of distinct subharmonic solutions (i.e. solutions with period $kT$ where $k \in \mathbb{N}$) for Hamiltonian systems with period $T$ have been explored in Rabinowitz \cite{Rabinowitz1980}, Tarentello \cite{Tarentello1988}, and Liu-Wang \cite{Liu1993}.

\section{Preservation of Symmetry}
\label{sec:preserv}
For the proof of preservation of symmetry, we first state a result from Mawhin, the proof of which is found in \cite{Mawhin1976}.
\begin{thm}\label{thm:MMT}
Given a real Hilbert space $H$, with inner product $\langle, \rangle$ and norm $|\cdot|$, let $L: \text{dom}(L)\subset H \rightarrow H$ be a linear, self-adjoint operator and $N: H \rightarrow H$ be a mapping with a bounded, linear G{\^ a}teaux derivative $N'$ on H such that, for each $x\in H$, $N'(x)$ is a symmetric operator. Denote by $\rho(A)$, $\sigma(A)$, and $r_\sigma(A)$ the resolvent set, the spectrum, and the spectral radius, respectively, of any linear operator $A$ in $H$.

Suppose there exist real numbers $\lambda < \mu$ such that
\begin{equation*}
    (\lambda,\mu) \subset \rho(L), \quad \lambda,\mu \in \sigma(L)
\end{equation*}
and real numbers $p,q$ with
\begin{equation*}
    \lambda < q \leq p < \mu
\end{equation*}
such that, for each $x\in H$,
\begin{equation*}
    qI \leq N'(x) \leq pI.
\end{equation*}
Then, $L - N$ is one-to-one,
\begin{equation*}
    (L-N)(\text{dom}(L)) = H,
\end{equation*}
and $(L-N)^{-1}$ is globally Lipschitzian. 
\end{thm}
We now present our result on preservation of symmetry, which is inspired by \cite{Willem1989}.
\begin{thm}
Let $g:\mathbb{R} \rightarrow \mathbb{R}$ be a given $C^1$ function. Let $s$ be any positive integer. Suppose there exists $i\geq 1$ such that
\begin{equation}\label{eq:condition}
    \lambda_i < \min_{t\in \mathbb{R}} g'(u(t)) < \max_{t\in \mathbb{R}} g'(u(t)) < \lambda_{i+1}
\end{equation}
then, for every $f\in L^2[0,2\pi]$ such that $f$ is $2\pi/s$-periodic, all solutions of $\eqref{eq:problem}$ are $2\pi/ks$-periodic for some $k\in \mathbb{N}$ (and so is also $2\pi/s$ periodic). 
\end{thm}

\textbf{Remark 1.}
Note that the last part of the statement of Theorem 2 suggests that if the hypotheses are satisfied, then all solutions will retain the same periodicity. However, this does not mean that the minimal period will be the same -- in fact, it could be this period divided by a positive integer.

\textbf{Remark 2.}
In light of using Theorem 1 to prove Theorem 2, it turns out that there exists a unique solution. Hence, Theorem 2 can be rephrased as saying that there exists a unique solution to (\ref{eq:problem}) that preserves the symmetry for every periodic function $f \in L^2[0,2\pi]$.

\begin{proof}
Let $P$ be the orthogonal projector on the space of $2\pi/s$-periodic functions. Define the self-adjoint operator $L:D(L) \subset L^2 \rightarrow L^2$ by
\begin{align*}
    D(L) &= H^1_{per}[0,2\pi] \cap H^2[0,2\pi]\\
    Lu &= -u''.
\end{align*}
For $u\in L^2$, we write $u = v+w$ for $v = Pu$ and $w = (I-P)u$. Then, \eqref{eq:problem} is equivalent to
\begin{align*}
    Lv &= -Pg(v+w) + f\\
    Lw &= -(I-P)g(v+w).
\end{align*}
Since $v$ is periodic with minimal period $2\pi/ks$ for some $k\in \mathbb{N}$, then $g(v)$ is also periodic with minimal period $2\pi/ks$. Note that if $w\equiv 0$, then
\begin{equation*}
    L(0) = -(I-P)g(v+0) = -(I-P)g(v) = 0.
\end{equation*}
Hence, $w = 0$ is a solution to $\eqref{eq:problem}$. In view of \eqref{eq:condition}, applying Theorem $\ref{thm:MMT}$, with $\lambda = \lambda_i, \mu = \lambda_{i+1}$, $q = \min_{t\in \mathbb{R}} g'(u(t))$, and  $p = \max_{t\in \mathbb{R}} g'(u(t))$, we see that for every $v\in \text{Range}(P)$, there is exactly one solution $w\in \text{Range}(I-P)$. In particular, $w = 0$ for every $v$. Hence, $u = v$. \qed
\end{proof}

\section{Breaking of Symmetry}
\label{sec:break}
For breaking of symmetry, we present a result inspired by \cite{CostaFang2019}. In fact, the following is a generalization of \cite{CostaFang2019} to arbitrary coprime positive integers.
\begin{thm}
Let $g:\mathbb{R} \rightarrow \mathbb{R}$ be a given $C^1$ function satisfying
\begin{gather}
    |g(t)| \leq C_g, \ \text{for all $t\in \mathbb{R}$}
\end{gather}
and
\begin{gather}
    G(t) \rightarrow -\infty \ \text{as} \ |t|\rightarrow\infty,
\end{gather}
where $G(t) = \int_0^t g(s)\ ds$. Suppose there exists $t_0,t_1\in \mathbb{R}$ such that
\begin{gather}\label{eigeninequality2}
    \lambda_{r-1} < g'(t_0) < \lambda_r < g'(t_1) < \lambda_{r+1}
\end{gather}
for some positive integer $r$ and let $s$ be another positive integer, $\gcd(r,s)=1$. Then there exists $\hat{f}:\mathbb{R}\rightarrow\mathbb{R}$, $\hat{f}\in V_s$ such that $\eqref{eq:problem}$ has at least one solution which is not in $V_s$. 
\end{thm}

\textbf{Remark 3.}
Let us understand the difference between Theorem 2 and Theorem 3. In the case of Theorem 2, if the derivative of the non-linearity $g$ is bounded between consecutive eigenvalues, each of which correspond to the periodic eigenfunctions, then the period of the solution(s) (if they exist) must retain at least the same period as the forcing term $f$. In contrast, if the derivative "crosses" an eigenvalue, then Theorem 3 states that one can find an appropriate $f$ where the periodicity of the solution(s) is \emph{not} preserved.

\begin{proof}
Define quadratic forms on $V_s^\perp$ by
\begin{gather*}
    Q_i(h) = \int_0^{2\pi} \left(|h'|^2 - g'(t_i)|h|^2\right)\ dt, \ i=0,1.
\end{gather*}
Note that since $r$ and $s$ are coprime, then $\sin(rt),\cos(rt)\in V_s^\perp$. Additionally, from \eqref{eigeninequality2}, $Q_i, i = 0,1,$ are non-degenerate.

Since $g'$ is continuous, then there exists $\delta_0 > 0$ such that $\lambda_{r-1} < g'(t) < \lambda_r$ for all $t\in [t_0-\delta_0,t_0+\delta_0]$. Consequently,
\begin{gather*}
    (r-1)^2 < g'(\delta_0\sin(st) + t_0) < r^2, \ \forall t\in \real.
\end{gather*}
Letting $u^*_0(t) \coloneqq \delta_0\sin(st) + t_0$, we define the quadratic form
\begin{gather*}
    \Tilde{Q}_0(h) = \int_0^{2\pi} (|h'|^2 - g'(u_0^*)|h|^2)\ dt, \ h\in V_s^\perp,
\end{gather*}
which is non-degenerate and has Morse index $m_0$. Similarly, define $\Tilde{Q}_1(h)$, which has Morse index $m_1$, using $u^*_1(t) \coloneqq \delta_1\sin(st) + t_1$ where $\delta_1 > 0$ is such that $\lambda_r < g'(t) < \lambda_{r+1}$ for all $t\in [t_1-\delta_1, t_1+\delta_1]$.

The rest of the proof follows as in \cite{CostaFang2019} with $s$ taking the role of $\hat{p}$ and $r$ taking the role of $p$. \qed
\end{proof}

\section{Acknowledgement}
We would like to thank Professor David Costa for his valuable suggestions for our manuscript and his inspiring mentorship. 
% \bibliographystyle{elsarticle-num-names} 
% \bibliography{cas-refs}





% \end{thebibliography}

\end{document}